\documentclass[12pt]{article}
\def\C{{\mathbf{C}}}
\def\N{{\mathbf{N}}}
\def\bC{{\mathbf{\overline{C}}}}
\def\const{{\mathrm{const}}}
\begin{document}
\title{Meromorphic functions of one complex variable.
A survey}
\author{A. Eremenko and J. K. Langley}
\date{}
\maketitle
\begin{abstract}
This is an appendix to the English translation
of the book by A. A. Goldberg and I. V. Ostrovskii,
Distribution of values of meromorphic functions,
Moscow, Nauka, 1970. An English translation of this book
is to be published soon by the AMS. In this appendix
we survey the results obtained on the topics of
the book after 1970.
\end{abstract}

The literature on meromorphic functions\footnote
{If the domain is not specified explicitly, we mean 
meromorphic
functions in $\C$.} is very large.
There is a comprehensive survey \cite{goldsurv1}
that contains everything that was reviewed on the topic in the
Soviet ``Referativnyi Zhurnal'' in 1953--1970,
and a later large survey \cite{goldsurvey2}.
More recent surveys
\cite{drasurv}, \cite{haymansurv} and \cite{e15}
are shorter and have narrower scope.

Some books on specific topics in
the theory of meromorphic functions published after 1970
are \cite{cherryye}, \cite{CL}, \cite{HK}, 
\cite{lang}, \cite{minru}, \cite{petrenko}, \cite{rubel}, \cite{lo}.
A survey of the fast developing
subject of iteration of meromorphic functions is \cite{Be}.

Here we give a short survey of some results which are closely related
to the problems considered in this book. Thus we do not
include many
important topics, like
geometric theory
of meromorphic functions,
iteration, composition, differential and
functional equations, normal families,
Borel and Julia directions,
uniqueness theorems,  and most
regrettably, holomorphic curves and quasiregular mappings. 
\vspace{.1in}

Chapter I.
\vspace{.1in}

Mokhon'ko \cite{mokhonko} proved the following generalization
of (6.29), Ch.~I. {\em Let $R(w,z)$ be a rational function
of $w$ whose coefficients
are meromorphic functions $h(z)$ satisfying $T(r,h)=O(\phi(r)),$
where $\phi$ is a fixed positive increasing function on $[0,\infty)$.
Then for every
meromorphic function $f$ we have}
$$T(r,R(f(z),z))=\deg_wR\, T(r,f)+O(\phi(r)).$$
This is a purely algebraic result;
its proof uses only properties (6.5), (6.6), ($6.8^\prime$)
and the property $T(r,f^2)=2T(r,f)$.

Another result in the same direction is due to Eremenko \cite{eremenkoDE}.
{\em Let $F(u,v,z)$ be a polynomial in $u$ and $v$
whose coefficients are meromorphic
functions $h$ of $z$
satisfying $T(r,h)=O(\phi(r)),$ where $\phi$
is a function as above, and assume that
$F$ is irreducible over the algebraic closure of the field of
meromorphic functions. If $f$ and $g$ are meromorphic functions
satisfying 
\begin{equation}
\label{5}
F(f(z),g(z),z)\equiv 0,
\end{equation}
then}
$$\deg_uF\, T(r,f)=(\deg_vF+o(1))T(r,g)+O(\phi(r)).$$
Both theorems have applications in the analytic theory of
differential equations.

Vojta \cite{vojta} noticed a formal analogy between the definition
of the Nevanlinna characteristic and the definition of {\em height}
in number theory. This analogy extends quite far, and it has been
a source of many interesting results and conjectures in the recent
years. For example, under Vojta's analogy, 
Jensen's formula corresponds to the fundamental
theorem of arithmetic, while
the second main theorem without the
ramification term corresponds
to the Thue--Siegel--Roth theorem
on Diophantine approximation.  The full second main Theorem
corresponds to the famous
unproved conjecture in number theory
which is known as the ``$abc$-conjecture''.
The analogy extends to
the multi-dimensional generalizations of Nevanlinna theory, where it has been
especially fruitful \cite{minru,lang2}.
As an example of the influence of Vojta's analogy on the one-dimensional
value distribution theory, we mention the recent precise
results
on the error term in the second main theorem (see comments to Chapter~III).

Analogies between function theory and number theory
were noticed before Vojta, see,
for example \cite{narasimhan}. 
Osgood was the first to
view the second main
theorem as an analogue of the Thue--Siegel Roth theorem
\cite{osgood0,osgood,osgood2}.
See also \cite{gundersen,waring}.
\vspace{.1in}

Formula (2.6), Ch.~I is the basis of the so-called ``Fourier method''
in the theory of entire and meromorphic functions, which was
developed in the work of Rubel \cite{rubel}, Taylor, Miles and others.
The Fourier method is considered
as a substitute for the Weierstrass
representation, which is more effective in certain cases.
One of the main results achieved with this technique is the
theorem of Miles \cite{miles} on the quotient representation
of meromorphic functions:
{\em Every meromorphic function $f$ can be written as $f=g_1/g_0$,
where $g_j$ are entire functions, possibly with common zeros,
satisfying
$$T(r,g_j)\leq AT(Br,f),\quad j=0,1,$$
where $A$ and $B$ are absolute constants (independent of $f$).}
Moreover, for every $B>1$ there exists $A$ such that the above statement
holds
for every $f$. Simple examples show that one cannot in general take
$B=1$. This result was improved by Khabibullin \cite{habibullin1}:
{\em Let $f$ be a meromorphic function and $\epsilon>0$ a non-increasing
convex function on $[0,\infty)$. Then there esists a representation
$f=g_1/g_0$ such that $g_j$ are entire functions, and}
$$\ln(|g_1(z)|+|g_0(z)|)\leq \frac{A_\epsilon}{\epsilon(|z|)}T((1+\epsilon(|z|))|z|,f)+B_\epsilon.$$
Khabibullin also extended Miles's theorem to meromorphic functions
in $\C^n$. A survey of his results is \cite{habibullin2}.

The Fourier method is also one of the main tools
for the study of functions
which satisfy various
restrictions on the arguments of $a$-points
(see comments to Chapter~VI). A recent book on the Fourier method is
\cite{kondrat}.
\vspace{.1in}

Chapter II. 
\vspace{.1in}

Various properties of increasing functions
on $[0,\infty)$ play an important role in the theory of meromorphic functions.
By 1970 it became clear that the lower order $\lambda$
is at least as important as
the order $\rho$.
The main development since 1970 was
``localizing'' the notion of
order.

Let $\Phi$ be an unbounded increasing function. A sequence
$r_k\to\infty$ is called a sequence of {\em P\'olya peaks} (of the first kind)
of order $\mu$
if the inequalities
$$\Phi(rr_k)\leq (1+\epsilon_k)r^\mu\Phi(r_k),
\quad \epsilon_k<r<\epsilon_k^{-1}$$
hold with some $\epsilon_k\to 0$. P\'olya peaks were
formally introduced by
Edrei \cite{polyapeaks1}, though they were used by
Edrei and Fuchs already in 1963.
P\'olya peaks of order $\mu$ exist for all
$\mu$ in certain interval $[\lambda^*,\rho^*]$, which contains
the interval $[\lambda,\rho]$. The endpoints $\lambda^*$ and $\rho^*$
of this interval are called the order and lower order
in the sense of
P\'olya, respectively. The following formulas were given in \cite{drasinshea}:
$$\rho^*=\sup\{ p\colon\limsup_{x,A\rightarrow\infty}g(Ax)/A^pg(x)=\infty\}$$ and
$$\lambda^*=\inf\{ p\colon\liminf_{x,A\rightarrow\infty}
g(Ax)/A^pg(x)=0\}.$$
Most contemporary results on functions
of
finite lower order use P\'olya peaks.
\vspace{.1in}

Chapter III.
\vspace{.1in}

Milloux' inequality (Ch.~III, Theorem 2.4)
has led to a rich vein of results developing the value
distribution properties of meromorphic functions and their
derivatives, in which a decisive
role has been played by the paper \cite{Hay1} of Hayman.

Perhaps the most striking of the many results from
\cite{Hay1} is
{\em Hayman's alternative} (Ch.~III, Theorem~2.6):
if a function $f$ meromorphic in the plane omits a finite value $a$,
and its $k$th derivative $f^{(k)}$, for some $k \geq 1$,
omits a finite non-zero value $b$,
then $f$ is constant.
Two principal questions arising in connection with Hayman's alternative
are:
(i) whether a version of
Hayman's main inequality (Ch.~III, (2.23)) holds with
$N(r, 1/(f-a))$ replaced by $\overline{N}(r, 1/(f-a))$;
(ii) whether $f^{(k)}$ can be replaced by a more general term,
such as a linear differential
polynomial
\begin{equation}
F = L[f] = f^{(k)} + a_{k-1} f^{(k-1)} + \ldots + a_0 f ,
\label{lindiffpoly}
\end{equation}
with suitable coefficients $a_j$ of small growth compared to $f$.
A positive answer to (i) was given by Chen \cite{Chen}. Question (ii) was
answered affirmatively in \cite{Laha},
although
there do exist exceptional functions $f$, which may be determined
from the $a_j$, for which $f$ and $F-1$ have no zeros. A unified
approach to the questions from
\cite{Chen,Laha}
may be found in \cite{BL}. It was shown further
in \cite{Chuang,hin1} that if $P$ is a non-constant
differential polynomial in $f$, all
of whose terms have degree at least $2$ in $f$ and its derivatives, then
a version of Hayman's inequality holds with $f^{(k)}$ replaced by $P$,
and with $\overline{N}$ counting functions.

Question (i) is related to the issue of whether $f- a$ and
$f^{(k)} - b$ ($k \geq 1, b \neq 0$) can both fail to have simple zeros,
in analogy with the sharp result that a nonconstant meromorphic
function cannot be completely branched over five distinct values.
It has recently been shown \cite{NPZ} using normal family methods
that if $f$ is transcendental and
meromorphic in the plane with only multiple zeros then $f'$ takes
every finite non-zero value infinitely often (see also
\cite{BL,WangFang}).

The obvious example $f(z) = e^z$ shows that a transcendental entire
function $f$ may have the property that $f$ and
all its derivatives omit $0$: thus
the condition $b \neq 0$ is necessary in Hayman's alternative.
However, Hayman showed in \cite{Hay1} (see Ch.~III, Theorem~2.7) that
if $f$ is an entire function such that $f f''$ has no zeros
then $f(z) = \exp( A z + B)$ with $A,B$ constants: this follows from
applying (Ch.~III, Theorem 2.6) to $f/f'$. Clunie \cite{Clunie} established
the corresponding result with
$f''$ replaced by a higher derivative $f^{(k)}$,
which on combination with Hayman's theorem on
$ff''$ significantly improved earlier results of
P\'olya and others
on the zeros of entire functions
and their derivatives \cite{csillag,Pol,Saxer,Saxer2}.
The {\em Tumura-Clunie method}, as developed
in \cite{Clunie}, shows that if
$\Psi_k [g]$ is a differential polynomial of form $\Psi_k [g] = g^k + P_{k-1}[g]
$,
where $P_{k-1}[g]$ is a polynomial in $g$ and its derivatives of
total degree at most $k-1$, and with coefficients of small growth compared
to $g$, and if $g$ has few poles and $\Psi_k[g]$ has few zeros, then
$\Psi_k[g]$ admits a simple factorisation. The application to zeros of
$ff^{(k)}$ then follows by writing $f^{(k)}/f$ in the form
$\Psi_k[f'/f]$. Variants and
generalisations of the Tumura-Clunie method, allowing
functions $g$ with unrestricted poles, appear in a large number of papers
including notably
\cite{MuesSteinmetz,Tohge5,ZhangLi}.

It was conjectured by Hayman in \cite{Hay1} that if $f$ is meromorphic
in the plane and $f f^{(k)}$ has no zeros, for some $k \geq 2$,
then $f(z) = \exp( Az + B)$ or $f(z) = (Az + B)^{-m}$, with $A, B$ constants
and $m \in \N$. This was proved by Frank \cite{frank} for $k \geq 3$
and by Langley \cite{langley1} for $k=2$ (see also
\cite{Bergff'',La10}), while simple examples show that
no such result holds for $k=1$. Generalisations include replacing
$ff^{(k)}$ by $fF$, where $F$ is given by (\ref{lindiffpoly})
with rational coefficients
\cite{Bru,FH,La9,Stei1}, and by $ff'' - \alpha f'^2 , \alpha \in \C$
\cite{Ber5,La10,Mues2}.
Closely linked to Milloux' inequality and
Hayman's alternative is the question of whether
$G= f^n f' $ must take every finite non-zero value,
when $f$ is non-constant meromorphic and $n \in \N$, the connection
being that $(n+1) G$ is the derivative of $f^{n+1}$, which has
only multiple zeros and poles. Hayman proved in \cite{Hay1} that this is
the case for $n \geq 2$ and $f$ entire, and for $n \geq 3$ when
$f$ has poles. Following incremental results by a number of authors
\cite{Clu0,Henne,Mues3}, the definitive theorem in this direction was
proved in \cite{BE}: if $f$ is a transcendental meromorphic
function in the plane and $m > k \geq 1$ then $(f^m)^{(k)}$ has
infinitely many zeros. Here the result is proved first for
finite order,
and the infinite order case is then deduced by applying
a renormalisation method from
normal families \cite{Pang1,Pang2,Zalcman2}. The closely related question of
whether $f' + f^n $ may omit a finite value, which in turn is related
to the Tumura-Clunie method, is resolved in
\cite{BE,Hay1,Mues3,Stei2}. A more recent conjecture \cite{WYY}
asserts that if $f$ is a transcendental meromorphic function then
$ff^{(k)} $ takes every finite, non-zero value infinitely often: this
is known to be
true for $k=1$ \cite{BE}, and for $k = 2$ and $f$ entire \cite{langley2}.

A further development from \cite{Hay1} leads to two conjectures which
remain open. Hayman observed in \cite{Hay1} that since the derivative
$f'$ of a transcendental meromorphic function $f$ has only multiple
poles, it follows that $f'$ has at most one finite Picard value.
It was subsequently conjectured by Mues
\cite{Mues1} that the Nevanlinna deficiencies
of $f^{(m)}$ satisfy, for $m \geq 1$,
\begin{equation}
\sum_{a \in \C} \delta (a, f^{(m)}) \leq 1.
\label{muesconj}
\end{equation}
This was proved by Mues \cite{Mues1} for
$m \geq 2$, provided
all poles of $f$ are simple. In the
general case the best known upper bound
for the sum in (\ref{muesconj}) appears to be $4/3$ \cite{ishizaki,YL}.

The Mues conjecture (\ref{muesconj}) would follow from a positive
resolution of the {\em Gol'dberg conjecture} that, for a transcendental
meromorphic function $f$ and $k \geq 2$,
\begin{equation}
\overline{N}(r, f) \leq N(r, 1/f^{(k)}) + o( T(r, f))
\label{goldconj}
\end{equation}
as $r \to \infty$, possibly outside an exceptional set. This is in turn linked
to a classical result of P\'olya
\cite{Pol} that if $f$ has at least two distinct poles
then $f^{(k)}$ has at least one zero, for all sufficiently large $k$.
When $f$
has poles of multiplicity at most $k-1$ the
inequality (\ref{goldconj}) follows from
a lemma of Frank and Weissenborn
\cite{FW} (see also \cite{Wang3}), so that in particular
if $f$ has only simple poles then (\ref{muesconj}) is true for every
positive $m$.
A related inequality
was proved in
\cite{Frankgraz,Steiw}: if
$F = L[f] $ as in
(\ref{lindiffpoly}), where the $a_j$ are small functions compared to $f$,
then either
$$
N(r, F) \leq N(r, 1/F) + 2 N(r, f) + o( T(r, f) )
$$
outside a set of finite measure or
$f$ is a rational function in solutions of the homogeneous equation
$L[w] = 0$. This method is
connected to Steinmetz' proof \cite{Stei3} of the second main
theorem for small functions discussed below.
 No further results in the direction
of (\ref{goldconj}) appear to be known, although it was proved by
Langley in \cite{Lanew}
that if $f$ is meromorphic of finite order in the plane
and $f^{(k)}$ has finitely many zeros, for some $k \geq 2$, then $f$
has finitely many poles.

Examples abound of meromorphic functions with infinitely many poles such
that the first derivative has no zeros,
but it was proved in \cite{ELR} (see also \cite{BE})
that if $f$ is transcendental
and meromorphic with
$
\limsup_{r \to \infty} T(r, f)/r = 0
$
then $f'$ has infinitely many zeros: the corresponding result with
$\limsup$ replaced by $\liminf$ may be found in \cite{hinchliffe}.

We remark that
the above results have all been stated for functions meromorphic in the plane.
Those which are proved only using properties of the Nevanlinna
characteristic admit in some cases generalisation to functions of
sufficiently rapid growth in a disk \cite{Hay1} or a half-plane
\cite{LeO}. Some related results for functions of slower growth
in the disc appear
in \cite{sheasons,sons}.

An old conjecture of Nevanlinna was that one can replace 
constants $a_k$ in the second main theorem by meromorphic functions
$a_k(z)$ with the property $T(r,a_k)=o(T(r,f))$. Such functions
$a_k$ are usually called ``small targets''. For the case of an entire
function $f$ such a generalization was obtained by Chuang Chi Tai in 1964.
Much later, the second main
theorem without the ramification term, was proved
for meromorphic functions by Osgood,
who used methods from number
theory \cite{osgood0,osgood,osgood2}.
A substantial simplification
of Osgood's proof was made by Steinmetz
\cite{Stei3}, who also used a beautiful idea
of Frank and Weissenborn \cite{FW}.
Osgood and Steinmetz proved that
$$\sum_{k=1}^q N(r,(f-a_k)^{-1})\geq (q-2+o(1))T(r,f),$$
outside of the usual exceptional set.
The proof in \cite{Stei3} is simple and elegant;
and uses only manipulations with Wronski determinants and
the classical lemma on the logarithmic derivative.
This makes it suitable for generalizations
to holomorphic curves
\cite{osgood,osgood2,minru}.
However, this version of the second main theorem
does not take ramification into account.
Simple examples like $f(z)=e^z+z$ where $\delta(\infty,f)=\delta(z,f)=1$
and $N_1(r,f)\sim T(r,f)$ show that one cannot include the term $N_1$.
However, the following
form of the second main theorem holds with small targets:
$$\sum_{k=1}^q\overline{N}(r,(f-a_k)^{-1})\geq (q-2+o(1))T(r,f),$$
where $\overline{N}(r,(f-a_k)^{-1})$ is the usual function counting zeros
of $f-a$ disregarding multiplicity.
This result was recently obtained by Yamanoi \cite{yamanoi2}.
In \cite{yamanoi1} he separately treats the case of rational functions $a_k$
when the proof is technically simpler.
Yamanoi's proof is very complicated,
and it will be hard to generalize to holomorphic curves.
Surprisingly,
it uses Ahlfors's theory of covering surfaces
(and also algebraic geometry,
moduli spaces of curves,
and combinatorics).
The idea to bring Ahlfors's theory to this
context has its origin in the work of
Sauer \cite{sauer} who obtained a partial result for rational
small targets. One application of Yamanoi's generalization
of the second main theorem is the following.
{\em Suppose that $f$ and $g$ are meromorphic functions in $\C$
satisfying a relation of the form (\ref{5}).
If the genus of the curve $F(u,v)=0$ is greater than $1$,
then $T(r,f)+T(r,g)=O(\phi(r))$.}
This was conjectured by Eremenko  in 1982, and the important special case
that $\phi(r)=\ln r$, 
that is $F$ is a polynomial in all three variables,
was proved by Zaidenberg in 1990.

Now we turn to the classic setting.
The estimate in the lemma on the logarithmic derivative was improved
by Gol'dberg and Grinstein \cite{goldberggrinstein}:
$$m(r,f'/f)\leq \ln^+\{ T(\rho,f)(1-(r/\rho))^{-1}\}+\const,$$
where the constant depends on $f$.
Vojta's analogy (see comments to Chapter~I) stimulated new interest
in refined estimates for the logarithmic derivative, as
well as for the error term
$$S(r,f)=\sum_{j=1}^q m(r,a_j,f)+N_1(r,f)-2T(r,f)$$
in the second main theorem. 
Miles \cite{mileslld} derived from Gol'dberg's and Grinstein's
estimate the following. {\em Let $\psi$ be a continuous non-decreasing
function such that
$$\int_1^\infty\frac{dt}{t\psi(t)}<\infty.$$
Then for every meromorphic function $f$ we have
$$m(r,f'/f)\leq \ln^+\psi(T(r,f))+O(1),$$ 
outside an exceptional set of finite logarithmic measure.}
The strongest results on the error term $S(r,f)$ belong
to Hinkkanen \cite{hinkkanen}, for example:
$$S(r,f)\leq \ln^+\psi(T(r,f))+O(1),$$ 
outside an exceptional set of finite logarithmic measure,
where $\phi$ is as before. If one replaces $\psi$
by $t\psi(t)$, then both results will hold outside an exceptional
set of finite measure. All 
presently known results on
the error terms in one-dimensional Nevanlinna theory
are collected in the book
\cite{cherryye}. 

No analog of the second main theorem holds without an
exceptional set of $r$'s. This can be seen from the 
result of Hayman \cite{hayman1}: {\em Let $\{ E_k\}$
be 
closed sets of zero logarithmic capacity, and $\phi$ and $\psi$
arbitrary unbounded increasing functions. Then there exist
an entire function $f$ and a sequence $r_k\to\infty$
such that
$N(r_k,a,f)\leq \phi(r_k)\ln r_k$ for all $a\in E_k,$
while $T(r_k,f)>\psi(r_k)$.}
\vspace{.1in}

Chapter IV.
\vspace{.1in}

The result of Hayman we just cited shows that the set
$E_V(f)$ of Valiron exceptional values, which always
has zero capacity, can contain any $F_\sigma$ set
of zero capacity, but a complete description of possible
sets $E_V(f)$
is not known. 

For meromorphic functions $f$ of finite order,
Hyllengren \cite{hyllengren}
obtained a very precise description of the sets $E_V(f)$.
Let us say that a set $E$ satisfies the $H$-condition if there
exist a sequence $(a_n)$ of complex numbers and $\eta>0$
such that every point of $E$ belongs to infinitely many discs
$\{ w:|w-a_n|<\exp(-\exp(\eta n ))\}$.
{\em For every meromorphic function $f$ of finite order,
and every $x\in (0,1)$,
the set $E_V(x,f)$ of those $a\in\C$ for which $\Delta(a,f)>x$
satisfies the $H$-condition, and
vise versa, for every set $E$ satisfying the $H$-condition,
there exist an entire function $f$ and a number $x\in(0,1)$ such
that $\Delta(a,f)>x$, for all $a\in E$.} 
Notice that the $H$-condition is much stronger than
the condition of zero capacity. 

\vspace{.1in}

The first example of an entire function of finite order whose
deficiency depends on the choice of the origin was constructed
by Miles \cite{miles3}. The order of this function
was very large. Then Gol'dberg, Eremenko and Sodin \cite{GES2}
constructed such examples with any given order
greater than $5$. (For entire
functions of order less than $3/2$,
deficiencies are independent of the choice of the origin.)
\vspace{.1in}

Chapter V.
\vspace{.1in}

Put $L(r,f)=\ln\mu(r,f)/\ln M(r,f)$,
where $\mu(r,f)=\min\{ |f(z):|z|=r\}$. Corollary~3 on p. 232
says that $\limsup_{r\to\infty}L(r,f)\geq-[2\lambda]$.
Hayman \cite{haymanmm} showed that the same holds with
$-2[\lambda]$ replaced by $-2.19\ln\rho$ when $\rho$ is large
enough. For functions of infinite order, he proved
$$\limsup_{r\to\infty}\frac{L(r,f)}{\ln\ln\ln M(r,f)}\geq -2.19.$$
He also constructed examples of entire functions of large
finite order for which $\limsup_{r\to\infty}L(r,f)<-1$.
Then Fryntov \cite{fryntov}, answering a question of
Hayman, constructed entire functions
of any given order $\rho>1$ with the same property. 
Drasin \cite{drasinf} constructed entire functions
of order one, maximal type, with the property $M(r,f)\mu(r,f)\to 0$.
This may be contrasted with
a remarkable theorem of Hayman \cite{haymanmm2} which 
says that {\em if $f$ is an entire function of order one and normal
type, and $M(r,f)\mu(r,f)$ is bounded, then $f(z)=c\exp(az)$
for some constants $c$ and $a$.} 

Thus for an entire function of order at least one,
$\mu(r,f)$ can decrease at a higher rate than that of
increase of $M(r,f)$. The situation changes dramatically
if we consider the rate of decrease of $|f(z)|$ on an
unbounded {\em connected} set.
Hayman and Kjellberg \cite{HKj} proved
that {\em for every entire function $f$
and every $K>1$ all components of the set 
$\{ z:\ln |f(z)|+K\ln M(|z|,f)<0\}$ are bounded.}

Theorem $1.3^\prime$ has been the subject of many deep generalizations.
First we mention the famous ``spread relation'' of Baernstein
\cite{baernstein1}
\cite{baernstein2}
conjectured by Edrei in \cite{edrei1}:
{\em If $f$ is a meromorphic function of lower order $\lambda$, 
then for every $\epsilon>0$ there are arbitrary large values of $r$
such that the set of arguments $\theta$ where $|f(re^{i\theta})|>1$
has measure at least}
\begin{equation}
\label{1}
\min\left\{\frac{4}{\lambda}\arcsin\sqrt{\frac{\delta(\infty,f)}{2}},\;
2\pi\right\}-\epsilon.
\end{equation}
Similar sharp estimates of the measure of the set where
$$\ln|f(re^{i\theta})|>\alpha T(r,f)$$ were given in \cite{AB}.

Fryntov, Rossi and Weitsman \cite{FRW,FRW2} proved that under the assumptions of
the spread conjecture, the set $|f(re^{i\theta})|>1$ must
contain
an arc of length (\ref{1}). 
See also \cite{baer0} for the sharp lower estimate of the length
of the arcs in the set
$\{\theta:\ln|f(re^{i\theta}|>\alpha \ln M(r,f)\}.$

Extremal functions for the spread relation and its generalizations
were studied extensively, \cite{baer,EF1,EF2,shea}.

The new methods introduced by Baernstein \cite{baernstein1},
\cite{baernstein2} are based on the use of subharmonic functions,
and especially, on a new type of maximal function, the so-called
``star-function'', which turned out to be very useful in solving
a wide variety of extremal problems of function theory.
An  account of Baernstein's star function and its main applications is 
contained in the monograph \cite{HK}.

One important application of the spread relation
is the sharp estimate of the sum of deficiencies
of a meromorphic function of lower order $\lambda\leq 1$ \cite{edrei}.
{\em If a meromorphic function $f$ of lower order $\lambda$ has at least two
deficient values, then}
$$\sum_{a\in\bC}\delta(a,f)\leq\left\{\begin{array}{ll}1-\cos\pi\lambda,\;&0<\lambda\leq1/2\\
2-\sin\pi\lambda,\;&1/2<\lambda\leq 1.
\end{array}\right.,$$
The sharp estimate of the sum of deficiencies of a meromorphic function
in terms of its order or lower order $\lambda$ is still not established
for $\lambda>1$. The conjectured extremal functions are
described in \cite{DW}.

The results of \S 2 show that neither Nevanlinna nor Borel exceptional
values need be asymptotic values. On the other hand,
Picard exceptional values are asymptotic. A natural question arises,
whether any condition of smallness of $N(r,a,f)$ in comparison
with $T(r,f)$ will imply that $a$ is an asymptotic value.
The basic result belongs to the intersection of
the papers \cite{lo1}, \cite{hayman2}, and
\cite{e2}. 
{\em Let $f$ be a meromorphic function of lower order $\lambda\leq\infty$.
If the order of $N(r,a,f)$ is strictly less than
$\min\{ 1/2,\lambda\}$ then $a$ is an asymptotic value.} 
Example 3 on p. 249 shows that this condition is sharp,
if only the lower order of $f$ is taken into account. 
In \cite{e2}, a weaker sufficient condition for
$a$ to be an asymptotic value is given, that uses both the order
and lower order of $f$.  Hayman \cite{hayman2}
gives the following refined condition: {\em if 
$$T(r,f)-\frac{1}{2}r^{1/2}\int_r^\infty t^{-3/2} N(t,a,f)dt\to\infty,$$
then $a$ is an asymptotic value.}
\vspace{.1in}

The problem on p. 285 of optimal estimation of $\kappa(f)$ for functions
of lower order greater than $1$ is still open,
even for entire functions. It has been solved only for entire
functions with zeros on a ray \cite{HWradial}.

The best estimates known at this time
for entire and meromorphic functions
with fixed $\lambda>1$
are contained in \cite{milesshea1,milesshea2}.
They are derived from the following sharp
inequality which is obtained by the Fourier method:
$$\limsup\frac{N(r,0)+N(r,\infty)}{m_2(r,f)}\geq \sup_{\lambda^*\leq\lambda\leq
\rho^*}
\sqrt{2}\frac{|\sin\pi\lambda|}{\pi\lambda}\left\{
1+\frac{\sin2\pi\lambda}{2\pi\lambda}\right\}^{-1/2},$$
where $m_2$ is the $L^2$-norm,
$$m_2^2(r,f)=\frac{1}{2\pi}\int_0^{2\pi}(\ln|f(re^{i\theta})|)^2d\theta,$$
and $\lambda^*$ and $\rho^*$ are the order and lower order in the sense
of P\'olya, \cite{drasinshea} (see also comments to Chapter II).
\vspace{.1in}

In 1929, F. Nevanlinna \cite{fnev} found that meromorphic functions
of finite order satisfying $N_1(r,f)\equiv 0$ have the following properties:

a) $2\rho$ is an integer, $2\rho\geq 2$, 

b) all deficient values are asymptotic, and

c) all deficiencies are rational numbers with denominators
at most $2\rho$, 
and their sum equals $2$.

It was natural to conjecture that one of the conditions
\begin{equation}
\label{2}
N_1(r,f)=o(T(r,f)),
\end{equation}
or
\begin{equation}
\label{3}
\sum_{ a\in\bC} \delta(a,f)=2
\end{equation}
implies the properties a), b) and c).
Notice, that by the second main theorem, (\ref{3}) implies (\ref{2})
for functions of finite order.

It turns out that a strong form of this conjecture holds:
\vspace{.1in}

{\sc Small Ramification Theorem}. {\em If $f$ is
a meromorphic function of finite lower order with the
property $(\ref{2})$ then a), b) and c) hold, and:
\begin{equation}
\label{reg}
T(r,f)=r^\rho\ell(r),
\end{equation}
where $\ell$ is a slowly varying function in the sense of Karamata.}
\vspace{.1in}

As a corollary we obtain that conditions (\ref{2}) and (\ref{3})
for functions of finite lower order are equivalent.

This result has a long history which begins with theorems of
Pfluger and Edrei and Fuchs establishing the case of
entire functions
(see Corollary 2 on p. 315). Weitsman \cite{weitsman7} proved that
(\ref{3}) implies that the number of deficient values is at most
$2\rho$. Then Drasin \cite{drasin7}, \cite{drasin3}
proved that for functions of finite order
(\ref{3}) implies a), b) and c) and the regularity
condition (\ref{reg}).
Eremenko 
proposed a new potential-theoretic
method (see, for example, \cite{e15})
which finally led to a proof of
a simpler proof of Drasin's theorem. 
The small ramification
theorem in its present form stated above is proved in \cite{e10}. 

These results show that besides the defect relation, there
is an additional restriction on defects
of functions of finite order: (\ref{3}) implies that the number
of deficient values is finite and all defects are rational.

There are other restrictions as well.
Weitsman's theorem \cite{weitsman} says
that for functions $f$ of finite lower order
\begin{equation}
\label{1/3}
\sum_{a\in\bC}\delta(a,f)^{1/3}<\infty.
\end{equation}
The story of this theorem is described in Comments to Chapter VII
(see p. 576).
Weitsman's proof can actually be modified to produce an upper
bound depending only on the lower order of $f$.

The second restriction concerns functions of
finite order having a defect equal to $1$.
Lewis and Wu proved that such functions satisfy
$$\sum_{a\in\bC}\delta(a,f)^{1/3-\alpha}<\infty,$$
with some absolute constant $\alpha>0$.
Their proof gives $\alpha=2^{-264}$, which is far from what is expected.
(Lewis and Wu state their result for entire functions
but their proof applies to all functions with $\delta(a,f)=1$ for some
$a$.)

Examples of entire functions of
finite order with the property $\delta(a_n,f)\geq c^n$
for some $c\in (0,1)$ are
constructed in \cite{e11}, but a large gap remains
between these examples and the result of Lewis and Wu. 
\vspace{.1in}

Recent research on value
distribution meromorphic functions of the form
\begin{equation}\label{pot}
\sum\frac{c_k}{z-z_k},\quad\sum\frac{|c_k|}{|z_k|}<\infty
\end{equation}
was mainly concentrated on the functions with $c_k>0$.
Such functions are (complex conjugate to) gradients of subharmonic
functions of genus zero with discrete mass.
The main conjecture is that every
function of the form (\ref{pot})
has zeros.
This was proved in \cite{ELR} under the additional
assumption that $\inf c_k>0$.
\vspace{.1in}

Chapter VI
\vspace{.1in}

Entire functions whose zeros lie on (or are close to) 
finitely many 
rays were intensively studied. Under certain conditions, one can estimate
$\delta(0,f)$ from below, as in Corollary 4 on p. 350.
The strongest results in this direction belong to 
Hellerstein and Shea \cite{HS} and Miles \cite{miles2}.
One of the results of \cite{HS} says that
$\delta(0,f)>B_q(\theta_1,...,\theta_n)$ {\em for all entire functions of
genus $q$ with zeros on the rays $\arg z\in\{\theta_1,...,\theta_q)$,
and $B_q\to 1$ when $q\to\infty$ while the rays remain fixed.}  
In the case of one ray, they obtained
$B_q(\theta)=1-(\pi^2e^{-1}+o(1))/\ln q,\; q\to\infty$.
For entire functions of infinite order with zeros on
a ray, Miles \cite{miles2} proved that
$N(r,0,f)/T(r,f)\to 0$ as $r\to\infty$ avoiding an 
exceptional set of zero logarithmic density. However,
it may happen that $\delta(0,f)=0$ for such functions,
as Miles shows by an example constructed in the same paper.

Hellerstein and Shea \cite{HS} also
considered meromorphic functions
of finite order whose zeros
$\{ z_n\}$ and poles $\{ w_n\}$ lie
in opposite sectors $|\arg z_n|\leq \eta$
and $|\arg w_n-\pi|\leq \eta$, where
$0\leq \eta<\pi/6$. For such functions, they obtained a sharp
estimate of $\kappa(f)$ (definition on p. 285) from above.

For entire functions with zeros on finitely many rays,
there are relations between the order and lower order (see p. 344).
There relations were further investigated in \cite{miles4,gleizer2} and
\cite{qiao}. 

Miles \cite{miles4} considers the class of meromorphic functions $f$
whose zeros belong to a finite union of rays $X$
and poles belong to a finite union of rays $Y$,
where $X\cap Y=\emptyset$, and such that the exponent of convergence
of the union of zeros and poles is a given number $q$.
He then produces a non-negative integer $p=p(q,X,Y)$ such that
$$ \lim_{r\to\infty}\frac{T(r,f)}{r^p}=\infty\quad\mbox{if}\quad p>0\quad
\mbox{and}\quad\lim_{r\to\infty}\frac{T(r,f)}{\ln r}=\infty\quad\mbox{if}
\quad p=0,$$
and these growth estimates are sharp in the considered class.
The integer $p$ depends in a subtle way on the arithmetical
properties
of the arguments of the rays $X$ and $Y$, and this integer is in general
hard
to compute or estimate.

Gleizer \cite{gleizer2} considers entire functions with zeros on $n$ rays.
If $n=1$ or $n=2$, we have
$\rho\leq [\lambda]+n$, where $[\;]$ is the integer part. This follows
from Theorem~1.1, Chapter VI. However, if $n=3$, then the
difference
$\rho-\lambda$ can be arbitrarily large. In this case,
Gleizer proved that $[\rho]\leq 3([\lambda]+1).$ For arbitrary $n$,
Qiao \cite{qiao} proved that $\rho\leq 4^{q-1}([\lambda]+1)$.

In \cite{gleizer1,gleizer3}, Gleizer extended Theorem~4.1 by taking into account
not only the order but the lower order in the sense of P\'olya.
He used Baernstein's star-function.
\vspace{.1in}

There has been
remarkable progress in the problems considered in \S 5.
The conjecture of P\'olya and Wiman stated on p. 417 was proved by
Hellerstein and Williamson \cite{HW1,HW2}.
{\em If $f$ is a real entire function
such that all zeros of $ff'f^{\prime\prime}$ are real then $f$ belongs to
the Laguerre--P\'olya class.} In \cite{HW3} the same authors
with Shen
classified
all entire functions (not necessarily real) with the property
that $ff'f^{\prime\prime}$ has only real zeros. The classification
of meromorphic functions with the property that all their derivatives
have only real zeros was achieved by Hinkkanen
\cite{hinkk1,hinkk2,hinkk3,hinkk4,hinkk5}.

Sheil-Small \cite{sheilsmall} proved a conjecture of Wiman (1911),
that every real entire function of finite order with the property
that $ff^{\prime\prime}$ has only real zeros belongs to
the Laguerre--P\'olya class. Bergweiler, Eremenko and Langley
\cite{BEL} extended Sheil-Small's result to functions
of infinite order.
Then Langley \cite{langley4} extended this result
to the derivatives of higher orders:
{\em If $f$ is a real entire function of infinite
order, with finitely many non-real
zeros, 
then $f^{(k)}$ has infinitely many non-real zeros for every $k\geq 2$.}
  
For real entire functions of finite order with finitely many non-real zeros,
that do not belong to the Laguerre--P\'olya class,
Bergweiler and Eremenko \cite{BE1} proved that the number of
non-real zeros of $f^{(k)}$ tends to infinity as $k\to\infty$.
Together with Langley's result, this confirms another conjecture
of P\'olya (1943).
\vspace{.1in}

Chapter VII
\vspace{.1in}

The inverse problem (as stated on p. 487)
was completely solved by Drasin \cite{drasin1}.
A simplified proof is given in \cite{drasin2}.
The general idea is the same as in Chapter VII:
quasiconformal surgery and a version of the theorem of
Belinskii and Teichm\"uller are involved.
However, unlike in Chapter VII, Drasin does not construct
the Riemann surface spread over the sphere explicitly
but uses a more flexible technique.

Theorem 8.1 in Chapter VII actually gives a complete
solution of the inverse problem for finitely many deficient
values in the class of meromorphic functions of finite order.
(This was not known in 1970 when the book was written.
That condition~3 of this Theorem 8.1
is necessary follows from
the small ramification theorem above).

On the narrow inverse problem in the class of meromorphic functions
of finite order with infinitely many deficiencies, there is the
following result \cite{e6}:
\vspace{.1in}

{\em Let $\{ a\}$ be an arbitrary infinite countable subset of $\bC$,
and $\{\delta_n\}$ positive numbers
satisfying the following conditions:

(i) $\delta_n<1$,

(ii) $\sum_n\delta_n<2$, and

(iii) $\sum_n\delta_n^{1/3}<\infty$.

Then there exists a meromorphic function $f$ of (large)
finite order
such that $\delta(a_n,f)=\delta_n$, and $f$ has no other
deficient values.}
\vspace{.1in}

The order of this function depends on the quantities
in the right hand sides of (i), (ii) and (iii).
\vspace{.1in}

Conditions (ii) and (iii) are 
necessary because of the small ramification theorem, and Weitsman's theorem
(see comments to Chapter V above).
Condition (i) cannot be removed because of the Lewis and Wu theorem
stated above, but it is not known what
the precise condition on $\delta_n$ is, if $\delta_1=1$.

The class of meromorphic functions
with finitely many critical and asymptotic
values which was used in Chapter VII to investigate the
inverse problem
is interesting independently of
this application. 
Let us call this class $S$. The first general
result on functions of this class
belongs to Teichm\"uller,
who proved that the second main theorem
becomes an asymptotic equality for functions of
this class:
$$\sum_{j=1}^qm(r,a_j,f)+N_1(r,f)=2T(r,f)+Q(r,f),$$
where $a_j$ are all critical and asymptotic values.

Langley \cite{langley3} found that the growth of a function
$f\in S$ cannot be arbitrary:
$$c(f):=\liminf_{r\to\infty}\frac{T(r,f)}{\ln^2r}>0.$$
This constant $c(f)$ can be arbitrarily small,
but in the case that $f$
has only three critical and asymptotic values, we have \cite{e14}
$c(f)\geq \sqrt{3}/(2\pi)$ and this is best possible.  
On the other hand, there are no restrictions from above
on the growth
of functions of class $S$ \cite{merenkov}.

Class $S$ plays an important role in
holomorphic dynamics
(iteration of entire and meromorphic functions),
see, for example,
the survey \cite{Be}.
In \cite{CER} an application of almost periodic
ends is given. In \cite{DW} the method of Chapter VII is
extended to a new class of Riemann surfaces 
which the authors call ``Lindel\"ofian ends''.
The corresponding functions have
infinitely many critical
values and thus do not belong to the class $S$.

\vspace{.1in}

Appendix.
\vspace{.1in}

Govorov's original proof of the Paley conjecture was
a byproduct of his research on the Riemann Boundary problem with
infinite index \cite{govorov}. His theorem was generalized to
meromorphic functions by Petrenko, to subharmonic functions
in ${{\mathbf{R}^n}}$ by Dahlberg, and to entire functions of several complex
variables by Khabubullin, see his survey \cite{habibullin2}.

Petrenko introduced the following quantity which he called
the ``deviation of $f$ from the point $a$''.
\begin{equation}
\label{beta}
\beta(a,f)=\liminf_{r\to\infty}\frac{\ln^+ M(r,(f-a)^{-1})}{T(r,f)}.
\end{equation}
This differs from the defect in one respect: the uniform norm 
of $\ln^+|f(re^{i\theta})-a|^{-1}$ stands in the numerator instead of the
$L^1$ norm. Petrenko's generalization of Govorov's theorem
proved in the Appendix can be restated
as:
\begin{equation}
\label{paley}
\beta(a,f)\leq \pi\lambda
\end{equation}
for all meromorphic functions of lower order $\lambda$ and all $a\in\bC$.
This was the starting point of a study of deviations $\beta(a,f)$ by Petrenko
and others. The results obtained before 1978 are
summarized in his book \cite{petrenko}. The main difference between the
theory of deviations and the theory of defects is the absence of
a first main theorem: there is no simple
relation between $\beta(a,f)$ and solutions of
the equation $f(z)=a$.

We only present a sample of the results.
By analogy with defects, one can expect that the
set of exceptional values in the sense of Petrenko
$$P(f):=\{ a\in\bC:\beta(a,f)>0\}$$
is small. This is indeed the case: {\em for every meromorphic function $f$,
the set $P(f)$ has zero logarithmic capacity; for functions of
finite lower order it is at most countable (but may have the power of
continuum for functions of infinite lower order).} The following
analog of the Defect relation for functions of finite lower order was
established by Marchenko and Shcherba \cite{MS}:
$$\sum_{a\in\bC}\beta(a,f)\leq\left\{\begin{array}{ll}2\pi\lambda,\;&\lambda\geq 1/2
,\\                           \pi\lambda\csc\pi\lambda,\;&\lambda <1/2.
\end{array}\right.$$
Moreover, an analog of Weitsman's theorem (see comments to Chapter V) holds:
{\em for functions $f$ of finite lower order we have}
$$\sum_{a\in\bC}\beta(a,f)^{1/2}<\infty,$$
and the exponent $1/2$ is best possible.
A version Baernstein's spread relation also holds with deviations
instead of deficiencies \cite{petrenko}. It is worth mentioning here that,
according to Baernstein \cite{baernstein2},
the idea of introducing the star
function that led to the proof of
the spread relation occurred under the influence
of Petrenko's proof of (\ref{paley}).

The inverse problem for deviations turned out to be simpler
than the inverse problem for deficiencies.
A complete solution for functions of finite order is given in \cite{e6}:
{\em For every at most countable set $\{ a_n\}$ of points
and every sequence of positive numbers $\beta_n$ satisfying
the condition $\sum\beta_n^{1/2}<\infty,$ there exists
a meromorphic function $f$ of finite order such that
$\beta(a_n,f)=\beta_n$ and $\beta(a,f)=0$ for $a\notin\{ a_n\}$.}

In general, there is no relation between the sets $E_N(f)$ and $P(f)$:
{\em for every pair $(A,B)$ of at most countable subsets of $\bC$,
there exists a meromorphic function $f$ of any given non-zero order
such that $E_{N}(f)=A$ and $P(f)=B$} \cite{GES1,GES2}. 
On the other hand, if $T(2r,f)=O(T(r,f))$ then  $P(f)=E_V(f)$
\cite{e5}.

Bergweiler and Bock \cite{BB} found an analog of (\ref{paley}) for 
functions of infinite lower order. The idea was to replace $T(r,f)$
in the denominator of (\ref{beta}) by $A(r,f)$. Notice that
if one uses the Ahlfors definition of $T(r,f)$ then
$A(r,f)=dT(r,f)/d\ln r$, for example, if $T(r,f)=r^\lambda$ then
$A(r,f)=\lambda T(r,f).$  Bergweiler and Bock
proved that {\em for every meromorphic function
$f$ of order at least $1/2$ and every $a\in\bC$ we have}
$$b(a,f):=\liminf_{r\to\infty}\frac{\ln^+M(r,(f-a)^{-1})}{A(r,f)}\leq \pi,$$
and then Eremenko \cite{e12} established the following analog
of the Defect Relation:
$$\sum_{a\in\bC}b(a,f)\leq\left\{\begin{array}{ll}2\pi,\;&\lambda>1/2,\\
                                    2\pi\sin\pi\lambda,\;&\lambda\leq 1/2
\end{array}\right.,$$
assuming that there are at least two values $a$ with $b(a,f)>0$.
It follows that for every meromorphic function the set $\{ a\in\bC:b(a,f)>0
\}$
is at most countable. 

Even Drasin's theorem on the extremal functions for the defect relation
(see Comments to Chapter V) has its analog for $b(a,f)$ \cite{e13}:

{\em If $f$ is of finite lower order and 
$$\sum_{a\in\bC}b(a,f)=2\pi$$ then
the following limit exists
$$\lim_{r\to\infty}\frac{\ln T(r,f)}{\ln r}=\frac{n}{2},$$
where $n$ is an integer, and $b(a,f)=\pi/n$ or $0$ for every $a\in \bC$.}

{\em A. E.: Purdue University

West Lafayette, IN 47907 USA
\vspace{.2in}

J. K. L.: University of Nottingham,

NG 7 2RD 

U. K.}
\end{document}